\documentclass{article}
\usepackage{amssymb,latexsym}
\usepackage{amsmath}
\usepackage{graphics}
\usepackage{graphicx}
\newtheorem{theorem}{Theorem}

\newtheorem{definition}{Definition}

\newtheorem{remark}{Remark}

\numberwithin{equation}{section} \numberwithin{theorem}{section}
\numberwithin{corollary}{section} \numberwithin{lemma}{section}
\numberwithin{remark}{section} \numberwithin{definition}{section}

\begin{document}
\title{Generalized Timelike Mannheim Curves in Minkowski space-time $E_1^4$}
\author{M. AKY\.{I}\v{G}\.{I}T, S. ERSOY, \.{I}. \"{O}ZG\"{U}R, M.
TOSUN}

\date{}

\maketitle
\begin{center}

Department of Mathematics, Faculty of Arts and Sciences \\
Sakarya University, 54187 Sakarya/TURKEY \end{center}

\begin{abstract}
We give a definition of generalized timelike Mannheim curve in
Minkowski space-time $E_1^4$. The necessary and sufficient
conditions for the generalized timelike Mannheim curve obtain. We
show some characterizations of generalized Mannheim curve.


\end{abstract}

\section{Introduction}\label{S:intro}

The geometry of curves has long captivated the interests of
mathematicians, from the ancient Greeks through to the era of
Isaac Newton (1647-1727) and the invention of the calculus. It is
branch of geometry that deals with smooth curves in the plane and
in the space by methods of differential and integral calculus. The
theory of curves is the simpler and narrower in scope because a
regular curve in a Euclidean space has no intrinsic geometry. One
of the most important tools used to analyze curve is the Frenet
frame, a moving frame that provides a coordinate system at each
point of curve that is "best adopted" to the curve near that
point. Every person of classical differential geometry meets early
in his course the subject of Bertrand curves, discovered in 1850
by J. Bertrand. A Bertrand curve is a curve such that its
principal normals are the principal normals of a second curve.
There are many works related with Bertrand curves in the Euclidean
space and Minkowski space, \cite{Bal1}-\cite{Mat1}.

\noindent Another kind of associated curve is called Mannheim
curve and Mannheim partner curve. The notion of Mannheim curves
was discovered by A. Mannheim in 1878. These curves in Euclidean
3-space are characterized in terms of the curvature and torsion as
follows: A space curve is a Mannheim curve if and only if its
curvature $\kappa$ and torsion $\tau$ satisfy the relation

\begin{equation*}
{k_1} = \beta \left( {k_1^2 + k_2^2} \right)
\end{equation*}
for some constant $\beta$. The articles concerning Mannheim curves
are rather few. In \cite{Blu}, a remarkable class of Mannheim
curves is studied. General Mannheim curves in the Euclidean
3-space are obtained in  \cite{Tig}-\cite{Orb} Recently, Mannheim
curves are generalized and some characterizations and examples of
generalized Mannheim curves are given in Euclidean 4-space $E^4$
by \cite{Mat2}.

\noindent In this paper, we study the generalized spacelike
Mannheim partner curves in  $4-$dimensional Minkowski space-time.
We will give the necessary and sufficient conditions for the
generalized spacelike Mannheim partner curves.

\section{Preliminaries}\label{S:intro}

To meet the requirements in the next sections, the basic elements
of the theory of curves in Minkowski space-time $E_1^4$ are
briefly presented in this section. A more complete elementary
treatment can be found in \cite{One}.

\noindent Minkowski space-time $E_1^4$ is an Euclidean space
provided with the standard flat metric given by

\begin{equation*}
\left\langle {\,\,,\,\,} \right\rangle  =  - dx_1^2 + dx_2^2 +
dx_3^2 + dx_4^2
\end{equation*}
where $\left( {{x_1},\,{x_2},\,{x_3},\,{x_4}} \right)$ is a
rectangular coordinate system in $E_1^4$.

\noindent Since $\left\langle {\;,\;} \right\rangle $ is an
indefinite metric, recall that a ${\bf{v}} \in E_1^4$ can have one
of the three causal characters; it can be spacelike if
$\left\langle {{\bf{v}},{\bf{v}}} \right\rangle  > 0$ or ${\bf{v}}
= 0$, timelike if $\left\langle {{\bf{v}},{\bf{v}}} \right\rangle
< 0$ and null(ligthlike) if $\left\langle {{\bf{v}},{\bf{v}}}
\right\rangle  = 0$ and ${\bf{v}} \ne 0$. Similarly, an arbitrary
curve ${\bf{c}} = {\bf{c}}\left( t \right)$ in $E_1^4$ can locally
be spacelike, timelike or null (lightlike) if all of its velocity
vectors ${\bf{c'}}\left( t \right)$ are, respectively, spacelike,
timelike or null. The norm of ${\bf{v}} \in E_1^4$ is given by
$\left\| {\bf{v}} \right\| = \sqrt {\left| {\left\langle
{{\bf{v}},{\bf{v}}} \right\rangle } \right|} $. If $\left\|
{{\bf{c'}}\left( t \right)} \right\| = \sqrt {\left| {\left\langle
{{\bf{c'}}\left( t \right),{\bf{c'}}\left( t \right)}
\right\rangle } \right|}  \ne 0$ for all $t \in I$, then $C$ is a
regular curve in $E_1^4$. A timelike (spacelike) regular curve $C$
is parameterized by arc-length parameter $t$ which is given by
${\bf{c}}:I \to E_1^4$, then the tangent vector ${\bf{c'}}\left( t
\right)$ along $C$ has unit length, that is, $\left\langle
{{\bf{c}}\left( t \right),{\bf{c}}\left( t \right)} \right\rangle
= 1\,,\,\,\,\left( {\left\langle {{\bf{c}}\left( t
\right),{\bf{c}}\left( t \right)} \right\rangle =  - 1} \right)$
for all $t \in I$.

\noindent Hereafter, curves considered are timelike and regular
${C^\infty }$ curves in $E_1^4$. Let $T\left( t \right) =
{\bf{c'}}\left( t \right)$ for all $t \in I$, then the vector
field $T\left( t \right)$ is timelike and it is called timelike
unit tangent vector field on $C$.

\noindent The timelike curve $C$ is called special timelike Frenet
curve if there exist three smooth functions
${k_1},\,\;{k_2},\;\,{k_3}$ on $C$ and smooth non-null frame field
$\left\{ {T,\,N,\,{B_1},\,{B_2}} \right\}$ along the curve $C$.
Also, the functions ${k_1},\,{k_2}$ and  ${k_3}$ are called the
first, the second and the third curvature function on $C$,
respectively. For the ${C^\infty }$ special timelike Frenet curve
$C$, the following Frenet formula is

\begin{equation*}
\left[ \begin{array}{l}
 {T'} \\
 {N'} \\
 {B_1}^\prime  \\
 {B_2}^\prime  \\
 \end{array} \right] = \,\left[ \begin{array}{l}
 \,\,\,0\,\,\,\,\,\,\,\,\,\,\,\,\,{k_1}\,\,\,\,\,\,\,\,\,\,\,\,\,0\,\,\,\,\,\,\,\,\,\,\,0 \\
 \,\,\,{k_1}\,\,\,\,\,\,\,\,\,\,\,0\,\,\,\,\,\,\,\,\,\,\,\,\,\,{k_2}\,\,\,\,\,\,\,\,\,0 \\
 \,\,\,0\,\,\,\,\,\,\,\, - {k_2}\,\,\,\,\,\,\,\,\,\,\,0\,\,\,\,\,\,\,\,\,\,{k_3} \\
 \,\,\,0\,\,\,\,\,\,\,\,\,\,\,\,\,\,0\,\,\,\,\,\,\,\, - {k_3}\,\,\,\,\,\,\,\,0 \\
 \end{array} \right]\left[ \begin{array}{l}
 T \\
 N \\
 {B_1} \\
 {B_2} \\
 \end{array} \right]
\end{equation*}
\cite{One}.

\noindent Here, due to characters of Frenet vectors of the
timelike curve, $T,\,N,\,{B_1}$ and ${B_2}$ are mutually
orthogonal vector fields satisfying equations
\begin{equation*}
\left\langle {T\,,\,T} \right\rangle  =  - 1\,\,,\,\,\left\langle
{N\,,\,N} \right\rangle  = \left\langle {{B_1}\,,\,{B_1}}
\right\rangle  = \left\langle {{B_2}\,,\,{B_2}} \right\rangle  =
1\,\,.\,\,
\end{equation*}
For $t \in I,$ the non-null frame field
$\left\{{T,\,N,\,{B_1},\,{B_2}} \right\}$ and curvature functions
${k_1},\,\;{k_2}$ and $\;{k_3}$ are determined as follows

\[\begin{array}{l}
 {1^{st}}\,\,\,\,\,\,{\rm{step   \,\,\, \,T}}\left( t \right) = {\bf{c'}}\left( t \right) \\
 {2^{nd}}\,\,\,\,{\rm{step   }}\,\,\,\,{k_1}\left( t \right) = \left\| {T'\left( t \right)} \right\| > 0 \\
 \,\,\,\,\,\,\,\,\,\,\,\,\,\,\,\,\,\,\,\,\,\,\,\,\,\,\,N\left( t \right) = \frac{1}{{{k_1}\left( t \right)}}T'\left( t \right) \\
 {3^{rd}}\,\,\,\,\,{\rm{step}}\,\,\,{k_2}\left( t \right) = \left\| {N'\left( t \right) - {k_1}\left( t \right)T\left( t \right)} \right\| > 0 \\
 \,\,\,\,\,\,\,\,\,\,\,\,\,\,\,\,\,\,\,\,\,\,\,\,\,\,{B_1}\left( t \right) = \frac{1}{{{k_2}\left( t \right)}}\left( {N'\left( t \right) - {k_1}\left( t \right)T\left( t \right)} \right) \\
 {4^{th}}\,\,\,\,\,{\rm{step}}\,\,\,{B_2}\left( t \right) = \varepsilon \frac{1}{{\left\| {{B_1}^\prime \left( t \right) + {k_2}\left( t \right)N\left( t \right)} \right\|}}\left( {{B_1}^\prime \left( t \right) + {k_2}\left( t \right)N\left( t \right)} \right) \\
 \end{array}\]

\noindent where $\varepsilon $ is determined by the fact that
orthonormal frame field $\left\{ {T\left( t \right),\,N\left( t
\right),} \right.$ $\left. {{B_1}\left( t \right),\,{B_2}\left( t
\right)} \right\}$, is of positive orientation. The function
${k_3}$ is determined by

\begin{equation*}
{k_3}\left( t \right) = \left\langle {{B_1}^\prime \left( t
\right)\,,\,{B_2}\left( t \right)} \right\rangle  \ne 0.
\end{equation*}

\noindent So the function ${k_3}$ never vanishes.

\noindent In order to make sure that the curve $C$ is a special
timelike Frenet curve, above steps must be checked, from
${1^{st}}$ step to ${4^{th}}$ step, for $t \in I.$

\noindent Let $\left\{ {T,\,N,\,{B_1},\,{B_2}} \right\}$ be the
moving Frenet frame along a unit speed timelike curve $C$ in
$E_1^4$, consisting of the tangent, the principal normal, the
first binormal and the second binormal vector field, respectively.
Since $C$ is a timelike curve, its Frenet frame contains only
non-null vector fields.

\section{Generalized timelike Mannheim curves in $E_1^4$}\label{S:intro}

Mannheim curves are generalized in   by \cite{Mat2}. In this
paper, we have investigated generalization of timelike Mannheim
curves in Minkowski space $E_1^4$.

\begin{definition}\label{T:3.1}
A special timelike curve $C$ in $E_1^4$ is a generalized timelike
Mannheim curve if there exists a special timelike Frenet curve
$C^*$ in $E_1^4$ such that the first normal line at each point of
$C$ is included in the plane generated by the second normal line
and the third normal line of ${C^*}$ at the corresponding point
under $\phi $. Here $\phi $ is a bijection from $C$ to ${C^*}$.
The curve ${C^*}$ is called the generalized timelike Mannheim mate
curve of $C$.
\end{definition}

\noindent By the definition, a generalized Mannheim mate curve
${C^*}$ is given by the map ${{\bf{c}}^*}:{I^*} \to E_1^4$ such
that

\begin{equation}
{{\bf{c}}^*}\left( t \right) = {\bf{c}}\left( t \right) + \beta
\left( t \right)N\left( t \right),\,\,t \in I.
\end{equation}
Here $\beta $ is a smooth function on $I.$ Generally, the
parameter $t$ isn't an arc-length of $C^*$. Let $t^*$ be the
arc-length of $C^*$ defined by

\begin{equation*}
{t^*} = \int\limits_0^t {\left\| {\frac{{d{{\bf{c}}^*}\left( t
\right)}}{{dt}}} \right\|dt.}
\end{equation*}
If a smooth function $f:I \to {I^*}$  is given by $f\left( t
\right) = {t^*}$, then for $\forall t \in I$, we have
\begin{equation*}
f'\left( t \right) = \frac{{d{t^*}}}{{dt}} = \left\|
{\frac{{d{{\bf{c}}^{\bf{*}}}\left( t \right)}}{{dt}}} \right\| =
\sqrt {\left| { - {{\left( {1 + \beta \left( t \right){k_1}\left(
t \right)} \right)}^2} + {{\left( {\beta '\left( t \right)}
\right)}^2} + {{\left( {\beta \left( t \right){k_2}\left( t
\right)} \right)}^2}} \right|} .
\end{equation*}
The representation of timelike curve ${C^*}$ with arc-length
parameter $t^*$ is
\begin{equation*}
\begin{array}{l}
 {{\bf{c}}^*}:\,{I^*} \to E_1^4 \\
 \,\,\,\,\,\,\,\,\,\,\,{t^*}\,\, \to \,{{\bf{c}}^*}\left( {{t^*}} \right). \\
 \end{array}
\end{equation*}
For a bijection $\phi :C \to {C^*}$ defined by $\phi \left(
{{\bf{c}}\left( t \right)} \right) = {{\bf{c}}^*}\left( {f\left( t
\right)} \right),$ the reparameterization of $C^*$ is

\begin{equation*}
{{\bf{c}}^*}\left( {f\left( t \right)} \right) = {\bf{c}}\left( t
\right) + \beta \left( t \right)N\left( t \right)
\end{equation*}
where $\beta $ is a smooth function on $I$. Thus, we have

\begin{equation*}
\frac{{d{c^*}\left( {f\left( t \right)} \right)}}{{dt}} = {\left.
{\frac{{d{c^*}\left( {{t^*}} \right)}}{{dt}}} \right|_{{t^*} =
f\left( t \right)}}f'\left( t \right) = f'\left( t
\right){T^*}\left( {f\left( t \right)} \right),\,\,t \in I.
\end{equation*}

\begin{theorem}\label{T:3.1}
If a special timelike Frenet curve $C$ in $E_1^4$ is a generalized
timelike Mannheim curve, then the following relation between the
first curvature function ${k_1}$ and the second curvature function
${k_2}$ holds:
\begin{equation}
{k_1}\left( t \right) =  - \beta \left( {k_1^2\left( t \right) -
k_2^2\left( t \right)} \right)\,\,,\,\,t \in I
\end{equation}
where $\beta $ is a constant number.
\end{theorem}

\noindent \textbf{Proof} Let $C$ be a generalized timelike
Mannheim curve and $C^*$ be the generalized timelike Mannheim mate
curve of $C$, as following diagram

\begin{equation*}
\begin{array}{l}
 \,\,\,\,\,\,\,\,\,\,\,\,\,\,\mathop {\bf{c}}\limits_{ \cdot \,\, \cdot } \,\,\,\,\,\,\,\,\,\,\,\,\,{\mathop {\bf{c}}\limits_{ \cdot \,\, \cdot } ^*}\, \\
 f:\,\,\,\,\,\,I\,\,\,\, \to \,\,\,{I^*}\, \\
 \,\,\,\,\,\,\,\,\,\,\,\,\,\,\, \downarrow \,\,\,\,\,\,\,\,\,\,\,\,\,\,\,\, \downarrow  \\
 \phi \,\,:\,\,\,E_1^4\,\, \to \,E_1^4 \\
 \end{array}
\end{equation*}
A smooth function $h$ is defined by $f\left( t \right) =
\int{\left\| {\frac{{d{{\bf{c}}^*}\left( t \right)}}{{dt}}}
\right\|} dt = {t^*}$ and ${t^*}$ is the arc-length parameter of
${C^*}.$ Also $\phi $ is a bijection which is defined by $\phi
\left({{\bf{c}}\left( t \right)} \right) = {{\bf{c}}^*}\left(
{f\left( t \right)} \right).$ Thus, the timelike curve ${C^*}$ is
reparametrized as follows

\begin{equation}
{{\bf{c}}^*}\left( {f\left( t \right)} \right) = {\bf{c}}\left( t
\right) + \beta \left( t \right)N\left( t \right)
\end{equation}
where $\beta :I \subset \mathbb{R} \to \mathbb{R}$ is a smooth
function. By differentiating both sides of equation (3.3) with
respect to $t$, we have
\begin{equation}
f'\left( t \right){T^*}\left( {f\left( t \right)} \right) = \left(
{1 + \beta \left( t \right){k_1}\left( t \right)} \right)T + \beta
'\left( t \right)N\left( t \right) + \beta \left( t
\right){k_2}\left( t \right){B_1}\left( t \right).
\end{equation}
On the other hand, since the first normal line at the each point
of $C$ is lying in the plane generated by the second normal line
and the third normal line of $C^*$ at the corresponding points
under bijection $\phi $, the vector field $N\left( t \right)$ is
given by
\begin{equation*}
N\left( t \right) = g\left( t \right){B_1}^*\left( {f\left( t
\right)} \right) + h\left( t \right){B_2}^*\left( {f\left( t
\right)} \right)
\end{equation*}
where $g$ and $h$  are some smooth functions on $I \subset
\mathbb{R}$. If we take into consideration
\begin{equation*}
\left\langle {{T^*}\left( {f\left( t \right)} \right),\,g\left( t
\right){B_1}^*\left( {f\left( t \right)} \right) + h\left( t
\right){B_2}^*\left( {f\left( t \right)} \right)} \right\rangle  =
0
\end{equation*}
and the equation (3.4), then we have $\beta '\left( t \right) = 0$
 So we rewrite the equation
(3.4) as
\begin{equation}
f'\left( t \right){T^*}\left( {f\left( t \right)} \right) = \left(
{1 + \beta {k_1}\left( t \right)} \right)T\left( t \right) + \beta
{k_2}\left( t \right){B_1}\left( t \right),
\end{equation}
that is,
\begin{equation*}
{T^*}\left( {f\left( t \right)} \right) = \frac{{\left( {1 + \beta
{k_1}\left( t \right)} \right)}}{{f'\left( t \right)}}T\left( t
\right) + \frac{{\beta {k_2}\left( t \right)}}{{f'\left( t
\right)}}{B_1}\left( t \right)
\end{equation*}
where
\begin{equation*}
f'\left( t \right) = \sqrt {\left| { - {{\left( {1 + \beta
{k_1}\left( t \right)} \right)}^2} + {{\left( {\beta {k_2}\left( t
\right)} \right)}^2}} \right|} .
\end{equation*}
By taking differentiation both sides of the equations (3.5) with
respect to $t \in I,$ we get

\begin{equation}
\begin{array}{l}
 f'\left( t \right)k_1^*\left( {f\left( t \right)} \right){N^*}\left( {f\left( t \right)} \right) = {\left( {\frac{{1 + \beta {k_1}\left( t \right)}}{{f'\left( t \right)}}} \right)^\prime }T\left( t \right) \\
 \,\,\,\,\,\,\,\,\,\,\,\,\,\,\,\,\,\,\,\,\,\,\,\,\,\,\,\,\,\,\,\,\,\,\,\,\,\,\,\,\,\,\,\,\,\,\,\, + \left( {\frac{{\left( {1 + \beta {k_1}\left( t \right)} \right){k_1}\left( t \right) - \beta {{\left( {{k_2}\left( t \right)} \right)}^2}}}{{f'\left( t \right)}}} \right)N\left( t \right) \\
 \,\,\,\,\,\,\,\,\,\,\,\,\,\,\,\,\,\,\,\,\,\,\,\,\,\,\,\,\,\,\,\,\,\,\,\,\,\,\,\,\,\,\,\,\,\,\,\,\, + {\left( {\frac{{\beta {k_2}\left( t \right)}}{{f'\left( t \right)}}} \right)^\prime }{B_1}\left( t \right) + \left( {\frac{{\beta {k_2}\left( t \right){k_3}\left( t \right)}}{{f'\left( t \right)}}} \right){B_2}\left( t \right). \\
 \end{array}
\end{equation}
Since
\begin{equation*}
\left\langle {{N^*}\left( {f\left( t \right)} \right),\,\,g\left(
t \right){B_1}^*\left( {f\left( t \right)} \right) + h\left( t
\right){B_2}^*\left( {f\left( t \right)} \right)} \right\rangle  =
0.
\end{equation*}
The coefficient of $N\left( t \right)$ in equation (3.6) vanishes,
that is,

\begin{equation*}
\left( {1 + \beta {k_1}\left( t \right)} \right){k_1}\left( t
\right) - \beta {\left( {{k_2}\left( t \right)} \right)^2} = 0.
\end{equation*}
Thus, this completes the proof.

\begin{theorem}\label{T:3.3}
In $E_1^4$, let $C$ be a special timelike Frenet curve such that
its non-constant first and second curvature functions satisfy the
equality ${k_1}\left( s \right) =  - \beta \left( {k_1^2\left( t
\right) - k_2^2\left( t \right)} \right)$ for all $t \in I \subset
\mathbb{R}.$ If the timelike curve ${C^*}$ given by
\begin{equation*}
{{\bf{c}}^*}\left( t \right) = {\bf{c}}\left( t \right) + \beta
N\left( t \right)
\end{equation*}
is a special timelike Frenet curve, then ${C^*}$ is a generalized
timelike Mannheim mate curve of $C$.
\end{theorem}
\noindent \textbf{Proof }The arc-length parameter of ${C^*}$ is
given by
\begin{equation*}
{t^*} = \int\limits_0^t {\left\| {\frac{{d{{\bf{c}}^*}\left( t
\right)}}{{dt}}} \right\|} dt\,\,,\,\,t \in I.
\end{equation*}
Under the assumption of
\begin{equation*}
{k_1}\left( t \right) =  - \beta \left( {k_1^2\left( t \right) -
k_2^2\left( t \right)} \right),
\end{equation*}
we obtain $f'\left( t \right) = \sqrt {\left| {1 + \beta
{k_1}\left( t \right)} \right|} \,,\,\,t \in I.$

\noindent Differentiating the equation ${{\bf{c}}^*}\left(
{f\left( t \right)} \right) = {\bf{c}}\left( t \right) + \beta
N\left( t \right)$ with respect to $t$ the we reach

\begin{equation*}
f'\left( t \right){T^*}\left( {f\left( t \right)} \right) = \left(
{1 + \beta {k_1}\left( t \right)} \right)T\left( t \right) + \beta
{k_2}\left( t \right){B_1}\left( t \right).
\end{equation*}
Thus, it is seen that
\begin{equation}
{T^*}\left( {f\left( t \right)} \right) = \left( {\frac{{1 + \beta
{k_1}\left( t \right)}}{{\sqrt {\left| {1 + \beta {k_1}\left( t
\right)} \right|} }}T\left( t \right) + \frac{{\beta {k_2}\left( t
\right)}}{{\sqrt {\left| {1 + \beta {k_1}\left( t \right)}
\right|} }}{B_1}\left( t \right)} \right),\,\,t \in I.
\end{equation}
The differentiation of the last equation with respect to $t$ is

\begin{equation}
\begin{array}{l}
 f'\left( t \right)k_1^*\left( {f\left( t \right)} \right){N^*}\left( {f\left( t \right)} \right) = {\left( {\sqrt {\left| {1 + \beta {k_1}\left( t \right)} \right|} } \right)^\prime }T\left( t \right) \\
 \,\,\,\,\,\,\,\,\,\,\,\,\,\,\,\,\,\,\,\,\,\,\,\,\,\,\,\,\,\,\,\,\,\,\,\,\,\,\,\,\,\,\,\, + \left( {\frac{{\left( {1 + \beta {k_1}\left( t \right)} \right){k_1}\left( t \right) - \beta k_2^2\left( t \right)}}{{\sqrt {\left| {1 + \beta {k_1}\left( t \right)} \right|} }}} \right)N\left( t \right) \\
 \,\,\,\,\,\,\,\,\,\,\,\,\,\,\,\,\,\,\,\,\,\,\,\,\,\,\,\,\,\,\,\,\,\,\,\,\,\,\,\,\,\,\,\,\, + {\left( {\frac{{\beta {k_2}\left( t \right)}}{{\sqrt {\left| {1 + \beta {k_1}\left( t \right)} \right|} }}} \right)^\prime }{B_1}\left( t \right) + \left( {\frac{{\beta {k_2}\left( t \right){k_3}\left( t \right)}}{{\sqrt {\left| {1 + \beta {k_1}\left( t \right)} \right|} }}} \right){B_2}\left( t \right). \\
 \end{array}
\end{equation}
From our assumption, we have

\begin{equation*}
\frac{{{k_1}\left( t \right) + \beta k_1^2\left( t \right) - \beta
k_2^2\left( t \right)}}{{\sqrt {\left| {1 + \beta {k_1}\left( t
\right)} \right|} }} = 0.
\end{equation*}
Thus, the coefficient of $N\left( t \right)$ in the equation (3.8)
is zero. It is seen from the equation (3.7), ${T^*}\left( {f\left(
t \right)} \right)$ is a linear combination of $T\left( t \right)$
and ${B_1}\left( t \right).$ Additionally, from equation (3.8),
${N^*}\left( {f\left( t \right)} \right)$ is given by linear
combination of $T\left( t \right),\;\,{B_1}\left( t \right)$ and
${B_2}\left( t \right)$. On the otherhand, ${C^*}$ is a special
timelike Frenet curve that the vector $N\left( t \right)$ is given
by linear combination of ${T^*}\left( {f\left( t \right)} \right)$
and ${N^*}\left( {f\left( t \right)} \right)$.

\noindent Therefore, the first normal line $C$ lies in the plane
generated by the second normal line and third normal line of $C^*$
at the corresponding points under a bijection $\phi $ which is
defined by $\phi \left( {{\bf{c}}\left( t \right)} \right) =
{{\bf{c}}^*}\left( {f\left( t \right)} \right)$.

\noindent This,completes the proof.

\begin{remark}\label{T:3.1}
In 4-diemsional Minkowski space $E_1^4$, a special timelike Frenet
curve $C$ with curvature functions  ${k_1}$ and ${k_2}$ satisfying
${k_1}\left( t \right) =  - \beta \left( {k_1^2\left( t \right) -
k_2^2\left( t \right)} \right)$, it is not clear that a smooth
timelike curve ${C^*}$ given by (3.1) is a special Frenet curve.
Thus, it is unknown whether the reverse of Theorem 3.1 is true or
false.
\end{remark}

\begin{theorem}\label{T:3.3}
Let $C$ be a special timelike  curve in $E_1^4$ with non-zero
third curvature function ${k_3}$. If there exists a timelike
special Frenet curve ${C^*}$ in $E_1^4$  such that the first
normal line of $C$  is linearly dependent with the third normal
line of $C^*$ at the corresponding points $c\left( t \right)$ and
${c^*}\left( t \right)$, respectively, under a bijection $\phi :C
\to {C^*}$, iff the curvatures ${k_1}$ and ${k_2}$ of $C$ are
constant functions.
\end{theorem}

\noindent \textbf{Proof} Let $C$ be a timelike Frenet curve in
$E_1^4$ with the Frenet frame field  $\left\{ {T,\,N,}
\right.$$\left. {{B_1},\,{B_2}} \right\}$ and curvature functions
${k_1},\,{k_2}$ and ${k_3}$. Also, we assume that $C^*$  be a
timelike special Frenet curve in $E_1^4$ with the Frenet frame
field $\left\{ {{T^*},\,{N^*},\,{B_1}^*,} \right.$ $\left.
{{B_2}^*} \right\}$ and curvature functions $k_1^*,\,\,k_2^*\,$
and $k_3^*$. Let the first normal line of $C$ be linearly
dependent with the third normal line of $C^*$ at the corresponding
points  $C$ and $C^*$, respectively. Then the parameterization of
$C^*$ is
\begin{equation}
{{\bf{c}}^*}\left( {f\left( t \right)} \right) = {\bf{c}}\left( t
\right) + \beta \left( t \right)N\left( t \right),\,\,t \in I.
\end{equation}
If the arc-length parameter of ${C^*}$ is given ${t^*}$, then

\begin{equation}
{t^*} = \int\limits_0^t {\sqrt {\left| { - {{\left( {1 + \beta
\left( t \right){k_1}\left( t \right)} \right)}^2} + \left( {\beta
'\left( t \right)} \right) + {{\left( {\beta \left( t
\right){k_2}\left( t \right)} \right)}^2}} \right|} } dt
\end{equation}
and
\begin{equation*}
\begin{array}{l}
 f:\,I \to {I^*} \\
 \,\,\,\,\,\,\,\,\,t\,\, \to \,f\left( t \right) = {t^*}. \\
 \end{array}
\end{equation*}
Moreover, $\phi :C \to {C^*}$ is a bijection given by $\phi \left(
{{\bf{c}}\left( t \right)} \right) = {{\bf{c}}^*}\left( {f\left( t
\right)} \right)$.

\noindent Differentiating the equation (3.9) with respect to $t$
and using Frenet formulas, we get
\begin{equation}
\begin{array}{l}
 f'\left( t \right){T^*}\left( {f\left( t \right)} \right) = \left( {1 + \beta \left( t \right){k_1}\left( t \right)} \right)T\left( t \right) + \beta '\left( t \right)N\left( t \right) \\
 \,\,\,\,\,\,\,\,\,\,\,\,\,\,\,\,\,\,\,\,\,\,\,\,\,\,\,\,\,\,\,\,\,\,\,\,\,\,\,\,\,\,\,\,\, + \beta \left( t \right){k_2}\left( t \right){B_1}\left( t \right). \\
 \end{array}
\end{equation}
Since  ${B_2}^*\left( {f\left( t \right)} \right) =  \mp N\left( t
\right)$, then
\begin{equation*}
\left\langle {f'\left( t \right){T^*}\left( {f\left( t \right)}
\right),\,{B_2}^*\left( {f\left( t \right)} \right)} \right\rangle
= \left\langle \begin{array}{l}
 \left( {1 + \beta \left( t \right){k_1}\left( t \right)} \right)T\left( t \right) + \beta '\left( t \right)N\left( t \right) \\
 \,\,\,\,\,\,\,\,\,\,\,\,\,\,\,\,\,\,\,\,\,\,\,\,\,\,\,\,\,\,\,\,\,\, + \beta \left( t \right){k_2}\left( t \right){B_1}\left( t \right),\, \mp N\left( t \right) \\
 \end{array} \right\rangle ,
\end{equation*}
that is,
\begin{equation*}
0 =  \mp \beta '\left( t \right).
\end{equation*}
From last equation, it is easily seen that $\beta $ is a constant.
Hereafter, we can denote $\beta \left( t \right) = \beta $, for
all $t \in I.$

\noindent From the equation (3.10), we have
\begin{equation*}
f'\left( t \right) = \sqrt {\left| { - {{\left( {1 + \beta
{k_1}\left( t \right)} \right)}^2} + {{\left( {\beta {k_2}\left( t
\right)} \right)}^2}} \right|}  > 0.
\end{equation*}
Thus, we rewrite the equation (3.11) as follows;
\begin{equation*}
{T^*}\left( {f\left( t \right)} \right) = \left( {\frac{{1 + \beta
{k_1}\left( t \right)}}{{f'\left( t \right)}}} \right)T\left( t
\right) + \left( {\frac{{\beta {k_2}\left( t \right)}}{{f'\left( t
\right)}}} \right){B_1}\left( t \right).
\end{equation*}
The differentiation of the last equation with respect to $t$ is

\begin{equation}
\begin{array}{l}
 f'\left( t \right)k_1^*\left( {f\left( t \right)} \right){N^*}\left( {f\left( t \right)} \right) = {\left( {\frac{{1 + \beta {k_1}\left( t \right)}}{{f'\left( t \right)}}} \right)^\prime }T\left( t \right) \\
 \,\,\,\,\,\,\,\,\,\,\,\,\,\,\,\,\,\,\,\,\,\,\,\,\,\,\,\,\,\,\,\,\,\,\,\,\,\,\,\,\,\,\,\,\,\,\,\,\,\,\,\,\,\,\,\, + \left( {\frac{{\left( {1 + \beta {k_1}\left( t \right)} \right){k_1}\left( t \right) - \beta k_2^2\left( t \right)}}{{f'\left( t \right)}}} \right)N\left( t \right) \\
 \,\,\,\,\,\,\,\,\,\,\,\,\,\,\,\,\,\,\,\,\,\,\,\,\,\,\,\,\,\,\,\,\,\,\,\,\,\,\,\,\,\,\,\,\,\,\,\,\,\,\,\,\,\,\,\, + {\left( {\frac{{\beta {k_2}\left( t \right)}}{{f'\left( t \right)}}} \right)^\prime }{B_1}\left( t \right) + \left( {\frac{{\beta {k_2}\left( t \right){k_3}\left( t \right)}}{{f'\left( t \right)}}} \right){B_2}\left( t \right). \\
 \end{array}
\end{equation}
Since $\left\langle {f'\left( t \right)k_1^*\left( {f\left( t
\right)} \right){N^*}\left( {f\left( t \right)}
\right),\,{B_2}^*\left( {f\left( t \right)} \right)} \right\rangle
= 0$ and  ${B_2}^*\left( {f\left( t \right)} \right) =  \mp
N\left( t \right)$ for all $t \in I$, we obtain

\begin{equation*}
{k_1}\left( t \right) + \beta k_1^2\left( t \right) - \beta
k_2^2\left( t \right) = 0
\end{equation*}
is satisfied. Then
\begin{equation}
\beta  =  - \frac{{{k_1}\left( t \right)}}{{k_1^2\left( t \right)
- k_2^2\left( t \right)}}
\end{equation}
is a non-zero constant number. Thus, from the equation (3.12), we
reach
\begin{equation*}
\begin{array}{l}
 {N^*}\left( {f\left( t \right)} \right) = \frac{1}{{f'\left( t \right)K\left( t \right)}}{\left( {\frac{{1 + \beta {k_1}\left( t \right)}}{{f'\left( t \right)}}} \right)^\prime }T\left( t \right) + \frac{1}{{f'\left( t \right)K\left( t \right)}}{\left( {\frac{{\beta {k_2}\left( t \right)}}{{f'\left( t \right)}}} \right)^\prime }{B_1}\left( t \right) \\
 \,\,\,\,\,\,\,\,\,\,\,\,\,\,\,\,\,\,\,\,\,\, + \frac{1}{{f'\left( t \right)K\left( t \right)}}\left( {\frac{{\beta {k_2}\left( t \right){k_3}\left( t \right)}}{{f'\left( t \right)}}} \right){B_2}\left( t \right) \\
 \end{array}
\end{equation*}
where $K\left( t \right) = k_1^*\left( {f\left( t \right)}
\right)$ for all $t \in I.$ Differentiating the last equation with
respect to $t$, then we have
\begin{equation*}
\begin{array}{l}
 f'\left( t \right)\left[ {k_1^*\left( {f\left( t \right)} \right){T^*}\left( {f\left( t \right)} \right) + k_2^*\left( {f\left( t \right)} \right){B_1}\left( {f\left( t \right)} \right)} \right] = {\left( {\frac{1}{{f'\left( t \right)K\left( t \right)}}{{\left( {\frac{{1 + \beta {k_1}\left( t \right)}}{{f'\left( t \right)}}} \right)}^\prime }} \right)^\prime }T\left( t \right) \\
 \,\,\,\,\,\,\,\,\,\,\,\,\,\,\,\,\,\,\,\,\,\,\,\,\,\,\,\,\,\,\,\,\,\,\,\,\,\,\,\,\,\,\,\,\,\,\,\,\, + \left( {\frac{{{k_1}\left( t \right)}}{{f'\left( t \right)K\left( t \right)}}{{\left( {\frac{{1 + \beta {k_1}\left( t \right)}}{{f'\left( t \right)}}} \right)}^\prime } - \frac{{{k_2}\left( t \right)}}{{f'\left( t \right)K\left( t \right)}}{{\left( {\frac{{\beta {k_2}\left( t \right)}}{{f'\left( t \right)}}} \right)}^\prime }} \right)N\left( t \right) \\
 \,\,\,\,\,\,\,\,\,\,\,\,\,\,\,\,\,\,\,\,\,\,\,\,\,\,\,\,\,\,\,\,\,\,\,\,\,\,\,\,\,\,\,\,\,\,\,\,\, + \left( {{{\left( {\frac{1}{{f'\left( t \right)K\left( t \right)}}{{\left( {\frac{{\beta {k_2}\left( t \right)}}{{f'\left( t \right)}}} \right)}^\prime }} \right)}^\prime } - \frac{{{k_3}\left( t \right)}}{{f'\left( t \right)K\left( t \right)}}\left( {\frac{{\beta {k_2}\left( t \right){k_3}\left( t \right)}}{{f'\left( t \right)}}} \right)} \right){B_1}\left( t \right) \\
 \,\,\,\,\,\,\,\,\,\,\,\,\,\,\,\,\,\,\,\,\,\,\,\,\,\,\,\,\,\,\,\,\,\,\,\,\,\,\,\,\,\,\,\,\,\,\,\,\, + \left( {{{\left( {\frac{1}{{f'\left( t \right)K\left( t \right)}}\left( {\frac{{\beta {k_2}\left( t \right){k_3}\left( t \right)}}{{f'\left( t \right)}}} \right)} \right)}^\prime } + \frac{{{k_3}\left( t \right)}}{{f'\left( t \right)K\left( t \right)}}{{\left( {\frac{{\beta {k_2}\left( t \right)}}{{f'\left( t \right)}}} \right)}^\prime }} \right){B_2}\left( t \right) \\
 \end{array}
\end{equation*}
for all  $t \in I.$ Considering

\begin{equation*}
\left\langle {f'\left( t \right)\left( {k_1^*\left( {f\left( t
\right)} \right){T^*}\left( {f\left( t \right)} \right) +
k_2^*\left( {f\left( t \right)} \right){B_1}^*\left( {f\left( t
\right)} \right)} \right)\,,\,\,{B_2}^*\left( {f\left( t \right)}
\right)} \right\rangle  = 0
\end{equation*}
and
\begin{equation*}
{B_2}^*\left( {f\left( t \right)} \right) =  \mp N\left( t
\right),
\end{equation*}
then we get

\begin{equation*}
{k_1}\left( t \right){\left( {\frac{{1 + \beta {k_1}\left( t
\right)}}{{f'\left( t \right)}}} \right)^\prime } - {k_2}\left( t
\right){\left( {\frac{{\beta {k_2}\left( t \right)}}{{f'\left( t
\right)}}} \right)^\prime } = 0.
\end{equation*}
Arranging the last equation, we find

\begin{equation}
\beta \left[ {{k_1}\left( t \right){k_1}^\prime \left( t \right) -
{k_2}\left( t \right){k_2}^\prime \left( t \right)}
\right]f'\left( t \right) - \left[ {{k_1}\left( t \right) + \beta
k_1^2\left( t \right) - \beta k_2^2\left( t \right)}
\right]f''\left( t \right) = 0.
\end{equation}
Moreover, the differentiation of the equation (3.13) with respect
to $t$ is
\begin{equation*}
{k'_1}\left( t \right) + 2\beta \left( {{k_1}\left( t
\right){{k'}_1}\left( t \right) - {k_2}\left( t
\right){{k'}_2}\left( t \right)} \right) = 0.
\end{equation*}
From the above equation, it is seen that
\begin{equation}
- \frac{{{{k'}_1}\left( t \right)}}{2} = \beta \left( {{k_1}\left(
t \right){{k'}_1}\left( t \right) - {k_2}\left( t
\right){{k'}_2}\left( t \right)} \right).
\end{equation}
Substituting the equations (3.13) and (3.15) into the equation
(3.14), we obtain
\begin{equation*}
 - \frac{{{{k'}_1}\left( t \right)}}{2} = 0.
\end{equation*}
This means that the first curvature function is constant (that is,
positive constant). Additionally, from the equation (3.15) it is
seen that the second curvature function ${k_2}$ is positive
constant, too.

\noindent Conversely, suppose that $C$ is a timelike Frenet curve
$E_1^4$ in with the Frenet frame field $\left\{
{T,\,N,\,{B_1},\,{B_2}} \right\}$ and curvature functions
${k_1},\,{k_2}$ and ${k_3}$. The first curvature function ${k_1}$
and the second curvature function ${k_2}$ of $C$ are of positive
constant. Thus, $\frac{{{k_1}}}{{k_2^2 - k_1^2}}$ is a positive
constant number, say $\beta.$

\noindent The representation of timelike curve  ${C^*}$ with
arc-length parameter $t$ is

\begin{equation}
\begin{array}{l}
 {{\bf{c}}^*}:\,I \to E_1^4 \\
 \,\,\,\,\,\,\,\,\,\,\,t\,\, \to \,\,{{\bf{c}}^*}\left( t \right) = {\bf{c}}\left( t \right) + \beta \left( t \right)N\left( t \right). \\
 \end{array}
\end{equation}
Let  ${t^*}$ denote the arc-length parameter of ${C^*}$, we have

\begin{equation*}
\begin{array}{l}
 f:\,I \to {I^*} \\
 \,\,\,\,\,\,\,t\,\, \to \,{t^*} = f\left( t \right) = \sqrt {\left| {1 + \beta {k_1}} \right|} t. \\
 \end{array}
\end{equation*}
Then, we obtain $f'\left( t \right) = \sqrt {\left| {1 + \beta
{k_1}} \right|}$ and

\begin{equation*}
\begin{array}{l}
 f'\left( t \right){T^*}\left( {f\left( t \right)} \right) = T\left( t \right) + \beta {N^\prime }\left( t \right) \\
 \,\,\,\,\,\,\,\,\,\,\,\,\,\,\,\,\,\,\,\,\,\,\,\,\,\,\,\,\,\,\,\,\,\, = \left( {1 + \beta {k_1}} \right)T\left( t \right) + \beta {k_2}{B_1}\left( t \right), \\
 \end{array}
\end{equation*}
that is

\begin{equation}
{T^*}\left( {f\left( t \right)} \right) = \sqrt {\left| {1 + \beta
{k_1}} \right|} T\left( t \right) + \frac{{\beta {k_2}}}{{\sqrt
{\left| {1 + \beta {k_1}} \right|} }}{B_1}\left( t \right).
\end{equation}
By differentiating both sides of the above equality with respect
to $t$ we find

\begin{equation*}
f'\left( t \right){\left. {\frac{{d{T^*}\left( {{t^*}}
\right)}}{{d{t^*}}}} \right|_{{t^*} = f\left( t \right)}} = \sqrt
{\left| {1 + \beta {k_1}} \right|} {T^\prime }\left( t \right) +
\frac{{\beta {k_2}}}{{\sqrt {\left| {1 + \beta {k_1}} \right|}
}}{B_1}^\prime \left( t \right)
\end{equation*}

\begin{equation*}
\begin{array}{l}
 \,\,\,\,\,\,\,\,\,\,\,\,\,\,\,\,\,\,\,\,\,\,\,\,\,\,\,\,\,\,\,\,\,\,\,\,\,\,\,\,\,\,\,\,\,\,\,\,\, = \left[ {\frac{{{k_1}\left( {1 + \beta {k_1}} \right) - \beta k_2^2}}{{\sqrt {\left| {1 + \beta {k_1}} \right|} }}} \right]N\left( t \right) + \left[ {\frac{{\beta {k_2}{k_3}\left( t \right)}}{{\sqrt {\left| {1 + \beta {k_1}} \right|} }}} \right]{B_2}\left( t \right) \\
  \\
 \,\,\,\,\,\,\,\,\,\,\,\,\,\,\,\,\,\,\,\,\,\,\,\,\,\,\,\,\,\,\,\,\,\,\,\,\,\,\,\,\,\,\,\,\,\,\,\,\, = \left[ {\frac{{\beta {k_2}{k_3}\left( t \right)}}{{\sqrt {\left| {1 + \beta {k_1}} \right|} }}} \right]{B_2}\left( t \right). \\
 \end{array}
\end{equation*}
Hence, since ${k_3}$ doesn't vanish, we get

\begin{equation*}
k_1^*\left( {f\left( t \right)} \right) = \left\| {{{\left.
{\frac{{d{T^*}\left( {{t^*}} \right)}}{{d{t^*}}}} \right|}_{{t^*}
= f\left( t \right)}}} \right\| = \varepsilon \frac{{\beta
{k_2}{k_3}\left( t \right)}}{{1 + \beta {k_1}}} > 0
\end{equation*}
where  $\varepsilon  = sign\left( {{k_3}} \right)$ denotes the
sign of function ${k_3}.$  That is, $\varepsilon $ is  $-1$ or
$+1$.

\noindent We can put

\begin{equation*}
{N^*}\left( {{t^*}} \right) = \frac{1}{{k_1^*\left( {{t^*}}
\right)}}\frac{{d{T^*}\left( {{t^*}} \right)}}{{d{t^*}}},\,\,t \in
I.
\end{equation*}
Then, we get

\begin{equation*}
{N^*}\left( {f\left( t \right)} \right) =  \mp {B_2}\left( t
\right).
\end{equation*}
Differentiating of the last equation with respect to $t$, we reach

\begin{equation*}
f'\left( t \right){\left. {\frac{{d{N^*}\left( {{t^*}}
\right)}}{{d{t^*}}}} \right|_{{t^*} = f\left( t \right)}} =  -
\varepsilon \frac{{{k_3}}}{{\sqrt {\left| {1 + \beta {k_1}}
\right|} }}{B_1}\left( t \right)
\end{equation*}
and we have

\begin{equation*}
f'\left( t \right){\left. {\frac{{d{N^*}\left( {{t^*}}
\right)}}{{d{t^*}}}} \right|_{{t^*} = f\left( t \right)}} -
k_1^*\left( {f\left( t \right)} \right){T^*}\left( {f\left( t
\right)} \right) =  - \varepsilon \frac{{\beta {k_2}{k_3}\left( t
\right)}}{{\sqrt {\left| {1 + \beta {k_1}} \right|} }}T\left( t
\right) - \varepsilon \sqrt {\left| {1 + \beta {k_1}} \right|}
{B_1}\left( t \right).
\end{equation*}
Since  $\varepsilon {k_3}\left( t \right)$ is positive for $t \in
I,$  we have

\begin{equation*}
\begin{array}{l}
 k_2^*\left( {f\left( t \right)} \right) = \left\| {{{\left. {\frac{{d{N^*}\left( {{t^*}} \right)}}{{d{t^*}}}} \right|}_{{t^*} = f\left( t \right)}} - k_1^*\left( {f\left( t \right)} \right){T^*}\left( {f\left( t \right)} \right)} \right\| \\
 \,\,\,\,\,\,\,\,\,\,\,\,\,\,\,\,\,\,\,\,\, = \,\,\sqrt { - \frac{{{\beta ^2}k_2^2{{\left( {{k_3}\left( t \right)} \right)}^2}}}{{1 + \beta {k_1}}} + \left( {1 + \beta {k_1}} \right){{\left( {{k_3}\left( t \right)} \right)}^2}}  \\
 \,\,\,\,\,\,\,\,\,\,\,\,\,\,\,\,\,\,\,\,\, = \,\,\sqrt {{{\left( {{k_3}\left( t \right)} \right)}^2}}  = \varepsilon {k_3}\left( t \right) > 0. \\
 \end{array}
\end{equation*}
Thus, we can put

\begin{equation*}
\begin{array}{l}
 {B_1}^*\left( {f\left( t \right)} \right) = \frac{1}{{k_2^*\left( {f\left( t \right)} \right)}}\left( {{{\left. {\frac{{d{N^*}\left( {{t^*}} \right)}}{{d{t^*}}}} \right|}_{{t^*} = f\left( t \right)}} - k_1^*\left( {f\left( t \right)} \right){T^*}\left( {f\left( t \right)} \right)} \right) \\
 \,\,\,\,\,\,\,\,\,\,\,\,\,\,\,\,\,\,\,\,\,\,\, = \,\, - \frac{{\beta {k_2}}}{{\sqrt {\left| {1 + \beta {k_1}} \right|} }}T\left( t \right) - \sqrt {\left| {1 + \beta {k_1}} \right|} {B_1}\left( t \right)\,\,,\,\,\,\,\,\,\,\,\,\,\,\,t \in I. \\
 \end{array}
\end{equation*}
Differentiation of the above with respect to $t$, we get
\begin{equation*}
f'\left( t \right){\left. {\frac{{d{B_1}^*\left( {{t^*}}
\right)}}{{d{t^*}}}} \right|_{{t^*} = f\left( t \right)}} =
\frac{{{k_2}}}{{\sqrt {\left| {1 + \beta {k_1}} \right|} }}N\left(
t \right) - {k_3}\left( t \right)\sqrt {\left| {1 + \beta {k_1}}
\right|} {B_2}\left( t \right).
\end{equation*}
Since $f'\left( t \right) = \sqrt {\left| {1 + \beta {k_1}}
\right|}$ and $k_2^*\left( {f\left( t \right)} \right){N^*}\left(
{f\left( t \right)} \right) = {k_3}\left( t \right){B_2}\left( t
\right)$, we have

\begin{equation*}
{\left. {\frac{{d{B_1}^*\left( {{t^*}} \right)}}{{d{t^*}}}}
\right|_{{t^*} = f\left( t \right)}} + k_2^*\left( {f\left( t
\right)} \right){N^*}\left( {f\left( t \right)} \right) =
\frac{{{k_2}}}{{1 + \beta {k_1}}}N\left( t \right).
\end{equation*}
Thus, we obtain ${B_2}^*\left( {f\left( t \right)} \right) =\delta
N\left( t \right)$ for $t \in I,$  where $\delta  =  \mp 1.$ We
must determine whether $\delta $ is $-1$ or $+1$ under the
condition that the frame field  $\left\{ {{T^*}\left( t
\right),\,\,{N^*}\left( t \right),\,\,{B_1}^*\left( t
\right),\,\,{B_2}^*\left( t \right)} \right\}$ is of positive
orientation.

\noindent We have, by $\det \left[ {T\left( t \right),\,N\left( t
\right),\,{B_1}\left( t \right),\,\,{B_2}\left( t \right)} \right]
= 1$ for $t \in I.$

\begin{equation*}
\begin{array}{l}
 \det \left[ {{T^*}\left( t \right),\,\,{N^*}\left( t \right),\,\,{B_1}^*\left( t \right),\,\,{B_2}^*\left( t \right)} \right] \\
 \,\,\,\,\, = \det \left[ \begin{array}{l}
 \sqrt {\left| {1 + \beta {k_1}} \right|} T\left( t \right) + \frac{{\beta {k_2}}}{{\sqrt {\left| {1 + \beta {k_1}} \right|} }}{B_1}\left( t \right), \\
 \,\,\,\,\,\,\,\,\,\,\,\,\,\varepsilon {B_2}\left( t \right), - \frac{{\beta {k_2}}}{{\sqrt {\left| {1 + \beta {k_1}} \right|} }}T\left( t \right) - \sqrt {\left| {1 + \beta {k_1}} \right|} {B_1}\left( t \right),\,\,\delta N\left( t \right)\, \\
 \end{array} \right] \\
 \,\,\,\,\, = \,\varepsilon \delta \left( {\left( {1 + \beta {k_1}} \right) - \frac{{{\beta ^2}k_2^2}}{{1 + \beta {k_1}}}} \right) = \varepsilon \delta  \\
 \end{array}
\end{equation*}
and  $\det \left[ {{T^*}\left( t \right),\,{N^*}\left( t
\right),\,{B_1}^*\left( t \right),\,\,{B_2}^*\left( t \right)}
\right] = 1$ for any $t \in I.$  Therefore,  we get  $\varepsilon
= \delta .$ Thus, we get

\begin{equation*}
{B_2}^*\left( {f\left( t \right)} \right) = \varepsilon N\left( t
\right)
\end{equation*}
and
\begin{equation*}
\begin{array}{l}
 k_3^*\left( {f\left( t \right)} \right) = \left\langle {{{\left. {\frac{{d{B_1}^*\left( {{t^*}} \right)}}{{d{t^*}}}} \right|}_{{t^*} = f\left( t \right)}},\,\,{B_2}^*\left( {f\left( t \right)} \right)} \right\rangle  \\
 \,\,\,\,\,\,\,\,\,\,\,\,\,\,\,\,\,\,\,\,\, = \,\,\varepsilon \frac{{{k_2}}}{{1 + \beta {k_1}}},\,\,\,\,\,t \in I. \\
 \end{array}
\end{equation*}
By the above facts,  ${C^*}$ is a special Frenet curve in $E_1^4$
and the first normal line at each point of  $C$ is the third
normal line of $C^*$ at corresponding each point under the
bijection $\phi :\,c \to \phi \left( {c\left( t \right)} \right) =
{c^*}\left( {f\left( t \right)} \right) \in {C^*}.$

\noindent Thus, the proof is completed.

\noindent The following theorem gives a parametric representation
of a generalized timelike Mannheim curves $E_1^4$.

\begin{theorem}
Let $C$ be a timelike special curve defined by

\begin{equation*}
{\bf{c}}\left( s \right) = \left[ {\begin{array}{*{20}{c}}
   {\beta \int {f\left( s \right)\cosh s\,\,ds} }  \\
   {\beta \int {f\left( s \right)\sinh s\,\,ds} }  \\
   {\beta \int {f\left( s \right)g\left( s \right)\,\,\,ds} \,}  \\
   {\beta \int {f\left( s \right)h\left( s \right)\,\,\,\,ds} \,\,}  \\
\end{array}} \right],\,\,s \in U \subset \mathbb{R}.
\end{equation*}
Here, $\beta $ is a non-zero constant number, $g:\,U \to
\mathbb{R}$ and $h:U \to \mathbb{R}$ are any smooth functions and
the positive valued smooth function $f:\,U \to \mathbb{R}$ is
given by

\begin{equation*}
\begin{array}{l}
 f = {\left( {1 - {g^2}\left( s \right) - {h^2}\left( s \right)} \right)^{{{ - 3} \mathord{\left/
 {\vphantom {{ - 3} 2}} \right.
 \kern-\nulldelimiterspace} 2}}}{\left( {1 - {g^2}\left( s \right) - {h^2}\left( s \right) + {{\dot g}^2}\left( s \right) + {{\dot h}^2}\left( s \right) - {{\left( {\dot g\left( s \right)h\left( s \right) - g\left( s \right)\dot h\left( s \right)} \right)}^2}} \right)^{{{ - 5} \mathord{\left/
 {\vphantom {{ - 5} 2}} \right.
 \kern-\nulldelimiterspace} 2}}} \\
 \,\,\,\,\,\,\,\,\,\,\,\,\,\,\,\,\,\,\,\,\,\,\,\,\left[ { - {{\left( {1 - {g^2}\left( s \right) - {h^2}\left( s \right) + {{\dot g}^2}\left( s \right) + {{\dot h}^2}\left( s \right) - {{\left( {\dot g\left( s \right)h\left( s \right) - g\left( s \right)\dot h\left( s \right)} \right)}^2}} \right)}^3}} \right. \\
 \quad \left. {\,\,\,\,\,\,\,\,\,\,\,\,\,\,\,\,\, + {{\left( {1 - {g^2}\left( s \right) - {h^2}\left( s \right)} \right)}^3}\left( \begin{array}{l}
  - {\left( {g\left( s \right) - \ddot g\left( s \right)} \right)^2} - {\left( {h\left( s \right) - \ddot h\left( s \right)} \right)^2} \\
  - {\left( {\left( {g\left( s \right)\dot h\left( s \right) - \dot g\left( s \right)h\left( s \right)} \right) - \left( {\dot g\left( s \right)\ddot h\left( s \right) - \ddot g\left( s \right)\dot h\left( s \right)} \right)} \right)^2} \\
  + {\left( {g\left( s \right)\ddot h\left( s \right) - \ddot g\left( s \right)h\left( s \right)} \right)^2} \\
 \end{array} \right)} \right], \\
 \end{array}
\end{equation*}
for $s \in U.$ Then the curvature functions ${k_1}$ and ${k_2}$ of
$C$ satisfy

\begin{equation*}
{k_1} =  - \beta \left( {k_1^2 - k_2^2} \right).
\end{equation*}
at the each point ${\bf{c}}\left( s \right)$ of  $C$.
\end{theorem}

\noindent \textbf{Proof }Let $C$  be a timelike special curve
defined by

\begin{equation*}
{\bf{c}}\left( s \right) = \left[ {\begin{array}{*{20}{c}}
   {\beta \int {f\left( s \right)\cosh s\,ds} }  \\
   {\beta \int {f\left( s \right)\sinh s\,ds} }  \\
   {\beta \int {f\left( s \right)g\left( s \right)\,\,ds} \,}  \\
   {\beta \int {f\left( s \right)h\left( s \right)\,\,ds\,} }  \\
\end{array}} \right]\quad ,\quad s \in U \subset \mathbb{R}
\end{equation*}
where  $\beta $ is a non-zero constant number, $g$ and $h$  are
any smooth functions. $f$  is a positive valued smooth function.
Thus, we obtain

\begin{equation}
\mathop c\limits^. \left( s \right) = \left[
{\begin{array}{*{20}{c}}
   {\beta f\left( s \right)\cosh s}  \\
   {\beta f\left( s \right)\sinh s}  \\
   {\beta f\left( s \right)g\left( s \right)}  \\
   {\beta f\left( s \right)h\left( s \right)}  \\
\end{array}} \right]\quad ,\quad s \in U \subset \mathbb{R}
\end{equation}
where the subscript prime $\left( . \right)$  denotes the
differentiation with respect to $s$.

\noindent The arc-length parameter $t$ of $C$  is given by

\begin{equation*}
t = \psi \left( s \right) = \int\limits_{{s_0}}^s {\left\|
{{\bf{\dot c}}\left( s \right)} \right\|} ds
\end{equation*}
where  $\left\| {{\bf{\dot c}}\left( s \right)} \right\| = \beta
f\left( s \right)\sqrt { - 1 + {g^2}\left( s \right) + {h^2}\left(
s \right)}$.

\noindent If $\varphi $ denotes the inverse function of $\psi :\,U
\to I \subset \mathbb{R}$, then $s = \varphi \left( t \right)$ and
we get

\begin{equation*}
\varphi '\left( t \right) = {\left\| {{{\left.
{\frac{{d{\bf{c}}\left( s \right)}}{{ds}}} \right|}_{s = \varphi
\left( t \right)}}} \right\|^{ - 1}}\quad ,\quad t \in I
\end{equation*}
where the prime $\left( ' \right)$ denotes the differentiation
with respect to $t$.

\noindent The unit tangent vector  $T\left( t \right)$ of the
curve $C$ at the each point ${\bf{c}}\left( {\varphi \left( t
\right)} \right)$
 is given by

\begin{equation}
T\left( t \right) = {\left( { - 1 + {g^2}\left( {\varphi \left( t
\right)} \right) + {h^2}\left( {\varphi \left( t \right)} \right)}
\right)^{ - {1 \mathord{\left/
 {\vphantom {1 2}} \right.
 \kern-\nulldelimiterspace} 2}}}\left[ {\begin{array}{*{20}{c}}
   {\cosh \left( {\varphi \left( t \right)} \right)}  \\
   {\sinh \left( {\varphi \left( t \right)} \right)}  \\
   {g\left( {\varphi \left( t \right)} \right)}  \\
   {h\left( {\varphi \left( t \right)} \right)}  \\
\end{array}} \right],\,\,t \in I.
\end{equation}
Some simplifying assumptions are made for the sake of brevity as
follows;

\begin{equation*}
\begin{array}{l}
 \sinh : = \sinh \left( {\varphi \left( t \right)} \right)\quad \quad ,\quad \,\,\cosh : = \cosh \left( {\varphi \left( t \right)} \right) \\
 f: = f\left( {\varphi \left( t \right)} \right)\quad \,\,\,\,\,\,\,\,\,\,\,\,\,\,\,,\quad \,\,g: = g\left( {\varphi \left( t \right)} \right)\quad \,\,\,\,\,\,\,\,\,\,\,\,,\quad h: = h\left( {\varphi \left( t \right)} \right), \\
 \dot g: = \dot g\left( {\varphi \left( t \right)} \right) = {\left. {\frac{{dg\left( s \right)}}{{ds}}} \right|_{s = \varphi \left( t \right)}}\quad \,\,\,\,\,\,\,\,\,\,\,\,\,\,\,\,,\quad \,\,\,\dot h: = \dot h\left( {\varphi \left( t \right)} \right) = {\left. {\frac{{dh\left( s \right)}}{{ds}}} \right|_{s = \varphi \left( t \right)}}, \\
 \ddot g: = \ddot g\left( {\varphi \left( t \right)} \right) = {\left. {\frac{{{d^2}g\left( s \right)}}{{d{s^2}}}} \right|_{s = \varphi \left( t \right)}}\quad \,\,\,\,\,\,\,\,\,\,\,\,\,,\quad \,\ddot h: = \ddot h\left( {\varphi \left( t \right)} \right) = {\left. {\frac{{{d^2}h\left( s \right)}}{{d{s^2}}}} \right|_{s = \varphi \left( t \right)}}, \\
 \varphi ': = \varphi '\left( t \right) = {\left. {\frac{{d\varphi }}{{dt}}} \right|_t}\,\,\,\,\,\,\,, \\
 A: = 1 - {g^2} - {h^2}\quad \,\,\,\,\,\,\,,\quad B: =  - g\dot g - h\dot h\quad \,\,\,\,\,\,\,\,\,\,,\quad C: =  - {{\dot g}^2} - {{\dot h}^2}, \\
 D: =  - g\ddot g - h\ddot h\quad \,\,\,\,\,\,\,\,\,\,\,,\quad E: =  - \dot g\ddot g - \dot h\ddot h\quad \,\,\,\,\,\,\,\,\,\,,\quad F: = {{\ddot g}^2} + {{\ddot h}^2}.\quad  \\
 \end{array}
\end{equation*}
Thus, we get

\begin{equation*}
\dot A = 2B\quad ,\quad \dot B = C + D\quad ,\quad \dot C =
2E\quad ,\quad \varphi ' = {\beta ^{ - 1}}{f^{ - 1}}{A^{{{ - 1}
\mathord{\left/
 {\vphantom {{ - 1} 2}} \right.
 \kern-\nulldelimiterspace} 2}}}.
\end{equation*}
So, we rewrite the equation (3.19) as

\begin{equation}
T: = T\left( t \right) = {A^{ - {1 \mathord{\left/
 {\vphantom {1 2}} \right.
 \kern-\nulldelimiterspace} 2}}}\left[ {\begin{array}{*{20}{c}}
   {\cosh }  \\
   {\sinh }  \\
   g  \\
   h  \\
\end{array}} \right].
\end{equation}
Differentiating the last equation with respect to $t,$  we find

\begin{equation*}
T' = \varphi '\left[ {\begin{array}{*{20}{c}}
   { - \frac{1}{2}{A^{ - {3 \mathord{\left/
 {\vphantom {3 2}} \right.
 \kern-\nulldelimiterspace} 2}}}\dot A\cosh  + {A^{ - {1 \mathord{\left/
 {\vphantom {1 2}} \right.
 \kern-\nulldelimiterspace} 2}}}\sinh }  \\
   { - \frac{1}{2}{A^{ - {3 \mathord{\left/
 {\vphantom {3 2}} \right.
 \kern-\nulldelimiterspace} 2}}}\dot A\sinh  + {A^{ - {1 \mathord{\left/
 {\vphantom {1 2}} \right.
 \kern-\nulldelimiterspace} 2}}}\cosh }  \\
   { - \frac{1}{2}{A^{ - {3 \mathord{\left/
 {\vphantom {3 2}} \right.
 \kern-\nulldelimiterspace} 2}}}\dot Ag + {A^{ - {1 \mathord{\left/
 {\vphantom {1 2}} \right.
 \kern-\nulldelimiterspace} 2}}}\dot g}  \\
   { - \frac{1}{2}{A^{ - {3 \mathord{\left/
 {\vphantom {3 2}} \right.
 \kern-\nulldelimiterspace} 2}}}\dot Ah + {A^{ - {1 \mathord{\left/
 {\vphantom {1 2}} \right.
 \kern-\nulldelimiterspace} 2}}}\dot h}  \\
\end{array}} \right],
\end{equation*}
that is,

\begin{equation}
T' =  - \varphi '{A^{ - {1 \mathord{\left/
 {\vphantom {1 2}} \right.
 \kern-\nulldelimiterspace} 2}}}\left[ {\begin{array}{*{20}{c}}
   {{A^{ - 1}}B\cosh  - \sinh }  \\
   {{A^{ - 1}}B\sinh  - \cosh }  \\
   {{A^{ - 1}}Bg - \dot g}  \\
   {{A^{ - 1}}Bh - \dot h}  \\
\end{array}} \right].
\end{equation}
From the last equation, we find

\begin{equation}
{k_1}: = {k_1}\left( t \right) = \left\| {T'\left( t \right)}
\right\| = \varphi '{A^{ - 1}}{\left( {A - AC + {B^2}} \right)^{{1
\mathord{\left/
 {\vphantom {1 2}} \right.
 \kern-\nulldelimiterspace} 2}}}.
\end{equation}
By the fact that $N\left( t \right) = {\left( {{k_1}\left( t
\right)} \right)^{ - 1}}T'\left( t \right)$, we get

\begin{equation*}
N: = N\left( t \right) =  - {A^{{1 \mathord{\left/
 {\vphantom {1 2}} \right.
 \kern-\nulldelimiterspace} 2}}}{\left( {A - AC + {B^2}} \right)^{{{ - 1} \mathord{\left/
 {\vphantom {{ - 1} 2}} \right.
 \kern-\nulldelimiterspace} 2}}}\left[ {\begin{array}{*{20}{c}}
   {{A^{ - 1}}B\cosh  - \sinh }  \\
   {{A^{ - 1}}B\sinh  - \cosh }  \\
   {{A^{ - 1}}Bg - \dot g}  \\
   {{A^{ - 1}}Bh - \dot h}  \\
\end{array}} \right].
\end{equation*}
In order to get second curvature function ${k_2}$, we need to
calculate ${k_2}\left( t \right) = \left\| {N'\left( t \right) -
{k_1}\left( t \right)T\left( t \right)} \right\|.$ After a long
process of calculations and using abbreviations, we obtain

\begin{equation}
N' - {k_1}T = \varphi '{A^{{{ - 3} \mathord{\left/
 {\vphantom {{ - 3} 2}} \right.
 \kern-\nulldelimiterspace} 2}}}{\left( {A - AC + {B^2}} \right)^{{{ - 3} \mathord{\left/
 {\vphantom {{ - 3} 2}} \right.
 \kern-\nulldelimiterspace} 2}}}\left[ {\begin{array}{*{20}{c}}
   {\left( {P + Q} \right)\cosh  - R\sinh }  \\
   {\left( {P + Q} \right)\sinh  - R\cosh }  \\
   {Pg - R\dot g + Q\ddot g}  \\
   {Ph - R\dot h + Q\ddot h}  \\
\end{array}} \right]
\end{equation}
where

\begin{equation}
\begin{array}{l}
 P = \left( {A - AC + {B^2}} \right)\left( {{B^2} - AC - AD} \right) - {\left( {A - AC + {B^2}} \right)^2} \\
 \,\,\,\,\,\,\,\,\,\, + AB\left( {B - AE + BD} \right), \\
 Q = {A^2}\left( {A - AC + {B^2}} \right), \\
 R = {A^2}\left( {B - AE + BD} \right). \\
 \end{array}
\end{equation}
If we simplify $P$ then we have

\begin{equation*}
P = {A^2}\left( {C - BE - D + CD - 1} \right).
\end{equation*}
Therefore, we rewrite the equations (3.23) and (3.24) as

\begin{equation}
N' - {k_1}T = \varphi '{A^{{{ - 1} \mathord{\left/
 {\vphantom {{ - 1} 2}} \right.
 \kern-\nulldelimiterspace} 2}}}{\left( {A - AC + {B^2}} \right)^{{{ - 3} \mathord{\left/
 {\vphantom {{ - 3} 2}} \right.
 \kern-\nulldelimiterspace} 2}}}\left[ {\begin{array}{*{20}{c}}
   {\left( {\tilde P + \tilde Q} \right)\cosh  - \tilde R\sinh }  \\
   {\left( {\tilde P + \tilde Q} \right)\sinh  - \tilde R\cosh }  \\
   {\tilde Pg - \tilde R\dot g + \tilde Q\ddot g}  \\
   {\tilde Ph - \tilde R\dot h + \tilde Q\ddot h}  \\
\end{array}} \right]
\end{equation}
where

\begin{equation}
\begin{array}{l}
 \tilde P = C - D + CD - BE - 1, \\
 \tilde Q = A - AC + {B^2}, \\
 \tilde R = B - AE + BD. \\
 \end{array}
\end{equation}
Consequently, from the equations (3.25) and (3.26), we have

\begin{equation*}
{\left\| {N' - {k_1}T} \right\|^2} = {\left( {\varphi '}
\right)^2}A{\left( {A - AC + {B^2}} \right)^{ - 3}}\left[
\begin{array}{l}
  - {\left( {\tilde P + \tilde Q} \right)^2} + {{\tilde R}^2} + {{\tilde P}^2}\left( {{g^2} + {h^2}} \right) + {{\tilde R}^2}\left( {{{\dot g}^2} + {{\dot h}^2}} \right) \\
  + {{\tilde Q}^2}\left( {{{\ddot g}^2} + {{\ddot h}^2}} \right) - 2\tilde P\tilde R\left( {g\dot g + h\dot h} \right) \\
  - 2\tilde R\tilde Q\left( {\dot g\ddot g + \dot h\ddot h} \right) + 2\tilde P\tilde Q\left( {g\ddot g + h\ddot h} \right) \\
 \end{array} \right].
\end{equation*}
Substituting the abbreviations into the last equation, we have

\begin{equation*}
{\left\| {N' - {k_1}T} \right\|^2} = {\left( {\varphi '}
\right)^2}A{\left( {A - AC + {B^2}} \right)^{ - 3}}[ - {\tilde
P^2}A - 2\tilde P\tilde Q - {\tilde Q^2} + {\tilde R^2} - {\tilde
R^2}C + {\tilde Q^2}F + 2\tilde P\tilde RB + 2\tilde R\tilde QE -
2\tilde P\tilde QD].
\end{equation*}
After substituting the equation (3.26) into the last equation and
simplifying it, we get

\begin{equation*}
\begin{array}{l}
 k_2^2 = {\left\| {N' - {k_1}T} \right\|^2} \\
 \quad  = {\left( {\varphi '} \right)^2}A{\left( {A - AC + {B^2}} \right)^{ - 2}}[\left( {A - AC + {B^2}} \right)\left( {1 - F} \right) + \left( {C - 1} \right){\left( {1 + D} \right)^2} - 2BE\left( {1 + D} \right) + A{E^2}]. \\
 \end{array}
\end{equation*}
Moreover, from the equation (3.22) it is seen that

\begin{equation*}
k_1^2 = {\left( {\varphi '} \right)^2}{A^{ - 2}}\left( {A - AC +
{B^2}} \right).
\end{equation*}
The last two equation gives us

\begin{equation*}
\begin{array}{l}
 k_2^2 - k_1^2 = {\left( {\varphi '} \right)^2}{A^{ - 2}}{\left( {A - AC + {B^2}} \right)^{ - 2}}\left[ { - {{\left( {A - AC + {B^2}} \right)}^3}} \right. \\
 \quad \left. { + {A^3}\left( {\left( {A - AC + {B^2}} \right)\left( {1 - F} \right) + \left( {C - 1} \right){{\left( {1 + D} \right)}^2} - 2BE\left( {1 + D} \right) + A{E^2}} \right)} \right]. \\
 \end{array}
\end{equation*}
By the fact $\varphi ' = {\beta ^{ - 1}}{f^{ - 1}}{A^{{{ - 1}
\mathord{\left/
 {\vphantom {{ - 1} 2}} \right.
 \kern-\nulldelimiterspace} 2}}}$, we obtain

\begin{equation}
\begin{array}{l}
 k_2^2 - k_1^2 = {\beta ^{ - 2}}{f^{ - 2}}{A^{ - 3}}{\left( {A - AC + {B^2}} \right)^{ - 2}}\left[ {{{\left( {A - AC + {B^2}} \right)}^3}} \right. \\
 \quad \quad \quad \;\left. { + {A^3}\left( {\left( {A - AC + {B^2}} \right)\left( {1 - F} \right) + \left( {C - 1} \right){{\left( {1 + D} \right)}^2} - 2BE\left( {1 + D} \right) + A{E^2}} \right)} \right] \\
 \end{array}
\end{equation}
and

\begin{equation*}
{k_1} = {\beta ^{ - 1}}{f^{ - 1}}{A^{ - {3 \mathord{\left/
 {\vphantom {3 2}} \right.
 \kern-\nulldelimiterspace} 2}}}{\left( {A - AC + {B^2}} \right)^{{1 \mathord{\left/
 {\vphantom {1 2}} \right.
 \kern-\nulldelimiterspace} 2}}}.
\end{equation*}
According to our assumption

\begin{equation*}
\begin{array}{l}
 f = {\left( {1 - {g^2} - {h^2}} \right)^{{{ - 3} \mathord{\left/
 {\vphantom {{ - 3} 2}} \right.
 \kern-\nulldelimiterspace} 2}}}{\left( {1 - {g^2} - {h^2} + {{\dot g}^2} + {{\dot h}^2} - {{\left( {\dot gh - g\dot h} \right)}^2}} \right)^{{{ - 5} \mathord{\left/
 {\vphantom {{ - 5} 2}} \right.
 \kern-\nulldelimiterspace} 2}}} \\
 \,\,\,\,\,\,\,\,\,\,\,\,\,\,\,\,\,\,\,\left[ { - {{\left( {1 - {g^2} - {h^2} + {{\dot g}^2} + {{\dot h}^2} - {{\left( {\dot gh - g\dot h} \right)}^2}} \right)}^3}} \right. \\
 \quad \left. {\,\,\,\,\,\,\,\,\,\,\,\,\,\,\, + {{\left( {1 - {g^2} - {h^2}} \right)}^3}\left( { - {{\left( {g - \ddot g} \right)}^2} - {{\left( {h - \ddot h} \right)}^2} - {{\left( {\left( {g\dot h - \dot gh} \right) - \left( {\dot g\ddot h - \ddot g\dot h} \right)} \right)}^2} + {{\left( {g\ddot h - \ddot gh} \right)}^2}} \right)} \right], \\
 \end{array}
\end{equation*}
we obtain

\begin{equation*}
f = {A^{ - {3 \mathord{\left/
 {\vphantom {3 2}} \right.
 \kern-\nulldelimiterspace} 2}}}{\left( {A - AC + {B^2}} \right)^{{{ - 5} \mathord{\left/
 {\vphantom {{ - 5} 2}} \right.
 \kern-\nulldelimiterspace} 2}}}\left[ \begin{array}{l}
 {\left( {A - AC + {B^2}} \right)^3} \\
 \,\,\,\,\,\,\,\, + {A^3}\left( \begin{array}{l}
 \left( {A - AC + {B^2}} \right)\left( {1 - F} \right) + \left( {C - 1} \right){\left( {1 + D} \right)^2} \\
  - 2BE\left( {1 + D} \right) + A{E^2} \\
 \end{array} \right) \\
 \end{array} \right].
\end{equation*}
Substituting the above equation into the equation into the
equations (3.27) and (3.28), we obtain

\begin{equation*}
{k_1} =  - \beta \left( {k_1^2 - k_2^2} \right).
\end{equation*}
The proof is completed.

\end {document}